\documentclass{article}
\usepackage{spconf,amsmath,graphicx}
\usepackage{amssymb}
\usepackage{epstopdf}
\usepackage{amsthm}
\usepackage{bm}
\usepackage{algorithm}
\usepackage[noend]{algpseudocode}
\usepackage{soul,color}
\usepackage{cite}
\usepackage{multirow}
\usepackage{fancyvrb} 
\VerbatimFootnotes    
\usepackage{listings}
\usepackage[normalem]{ulem}

\def\P{{\cal P}}

\def\E{{\cal E}}
\def\F{{\cal F}}

\def\R{{\mathbb R}}
\def\N{{\cal N}}

\newcommand{\hlnew}[1]{{#1}}

\newtheorem{definition}{Definition}
\newtheorem{remark}{Remark}
\newtheorem{theorem}{Theorem}
\newtheorem{lemma}{Lemma}

\newcommand{\norm}[1]{\|#1\|} 
\newcommand\abs[1]{\left|#1\right|}
\newcommand\numberthis{\addtocounter{equation}{1}\tag{\theequation}}

\DeclareMathOperator*{\vect}{vec}

\makeatletter
\newcommand\fs@betterruled{%
 \def\@fs@cfont{\bfseries}\let\@fs@capt\floatc@ruled
 \def\@fs@pre{\vspace*{5pt}\hrule height.8pt depth0pt \kern2pt}%
 \def\@fs@post{\kern2pt\hrule\relax}%
 \def\@fs@mid{\kern2pt\hrule\kern2pt}%
 \let\@fs@iftopcapt\iftrue}
\floatstyle{betterruled}
\restylefloat{algorithm}
\makeatother

\newcommand{\topnew}{\top \hspace{-0.05cm}}

\title{Exact Linear Convergence Rate Analysis for Low-Rank \hlnew{Symmetric} Matrix Completion via Gradient Descent}
\name{Trung Vu and Raviv Raich}
\address{School of EECS, Oregon State University, Corvallis, OR 97331-5501, USA \\ \{vutru, raich\}@oregonstate.edu}

\begin{document}
\ninept
\maketitle
\begin{abstract}
Factorization-based gradient descent is a scalable and efficient algorithm for solving low-rank matrix completion. Recent progress in structured non-convex optimization has offered global convergence guarantees for gradient descent under certain statistical assumptions on the low-rank matrix and the sampling set. However, while the theory suggests gradient descent enjoys fast linear convergence to a global solution of the problem, the universal nature of the bounding technique prevents it from obtaining an accurate estimate of the rate of convergence. In this paper, we perform a local analysis of the exact linear convergence rate of gradient descent for factorization-based matrix completion for symmetric matrices. Without any additional assumptions on the underlying model, we identify the deterministic condition for local convergence of gradient descent, which only depends on the solution matrix and the sampling set. More crucially, our analysis provides a closed-form expression of the asymptotic rate of convergence that matches exactly with the linear convergence observed in practice. To the best of our knowledge, our result is the first one that offers the exact rate of convergence of gradient descent for matrix factorization in Euclidean space for matrix completion.
\end{abstract}
%

\section{Introduction}
\label{sec:intro}

Matrix completion is the problem of recovering a low-rank matrix from a sampling of its entries. In machine learning and signal processing, it has a wide range of applications including collaborative filtering \cite{rennie2005fast}, system identification \cite{liu2009interior} and dimension reduction \cite{candes2011robust}. In the era of big data, matrix completion has been proven to be an efficient and powerful framework to handle the enormous amount of information by exploiting low-rank structure of the data matrix. 


Let $\bm M \in \R^{n \times m}$ be a rank $r$ matrix with \hlnew{$1 \leq r \leq \min(n,m)$}, and $\Omega = \{ (i,j) \mid M_{ij} \text{ is observed} \}$ be an index subset of cardinality $s$ such that $s \leq nm$. The goal is to recover the unknown entries of $\bm M$. Matrix completion can formulated as a linearly constrained rank minimization or a rank-constrained least squares problem \cite{candes2009exact}. Two popular approaches for solving the aforementioned matrix completion problem formulations are convex relaxation via nuclear norm and non-convex factorization. The former approach, motivated by the \hlnew{success} of compressed sensing, replaces the matrix rank by its convex surrogate (the nuclear norm). Extensive work on designing convex optimization algorithms with guarantees can be found in \cite{candes2009exact,ji2009accelerated,toh2010accelerated,cai2010singular,ma2011fixed}. 
While on the theoretical side, the solutions of the relaxed problem and the original problem can be shown to coincide with high probability, on the practical side, computational complexity concerns \hlnew{diminish} the applicability of these algorithms. When the size of the matrix grows rapidly, storing and optimizing over a matrix variable become computationally expensive and even infeasible. In addition, it is evident this approach suffers from slow convergence \cite{koren2009bellkor,vu2019accelerating}. 
In the second approach, the original rank-constrained optimization is studied. Interestingly, by reparametrizing the $n \times m$ matrix as the product of two smaller matrices $\bm M = \bm X \bm Y^\topnew$, for $\bm X \in \R^{n \times r}$ and $\bm Y \in \R^{m \times r}$, the resulting equivalent problem is unconstrained and more computationally efficient to solve \cite{burer2005local}. While this problem is non-convex, recent progress shows that for such problem 
any local minimum is also a global \hlnew{minimum} \cite{sun2016guaranteed,ge2016matrix}. Thus, basic optimization algorithms such as gradient descent \cite{chen2015fast,sun2016guaranteed,ma2018implicit} and alternating minimization \cite{chen2012matrix,jain2013low,hardt2014understanding,hardt2014fast} can provably solve matrix completion under a specific sampling regime. 
\hlnew{Alternatively, the original rank-constrained optimization problem can be solved without the aforementioned reparameterization via the truncated singular value decomposition.} \cite{jain2010guaranteed,goldfarb2011convergence,tanner2013normalized,jain2015fast,chunikhina2014performance,vu2019accelerating,vu2019local}. 

Among the aforementioned algorithms, let us draw our attention to the gradient descent method due to its outstanding simplicity and scalability. The first global convergence guarantee is attributed to Sun and Luo \cite{sun2016guaranteed}. The authors proved that gradient descent with appropriate regularization can converge to the global optima of a factorization-based formulation at a linear rate.
Later on, Ma~\textit{et.~al.} \cite{ma2018implicit} proposed that even in the absence of explicit regularization, gradient descent recovers the underlying low-rank matrix by implicitly regularizing its iterates. \hlnew{The} aforementioned results, while \hlnew{establishing} powerful guarantees on the convergence behavior of gradient descent, impose several limitations. 
For some methods, the linear convergence rate depends on \hlnew{constants that are not in closed-form and are} hard to verify in numerical \hlnew{experiments} even when the underlying matrix is known. Second, a solution-independent analysis of the convergence rate typically offers a loose bound when considered for a specific solution. 
Third, the global convergence analysis requires certain assumptions on the underlying model which largely restrict the setting of the matrix completion problem in practice. Among such assumptions, one would consider the \hlnew{incoherence} of the target matrix, the randomness of the sampling set, and the fact that the rank $r$ and the condition number of $\bm M$ are small constants as $n,m \to \infty$.

To address these issues, we consider the local convergence analysis of gradient descent for factorization-based matrix completion. \hlnew{In the scope of this paper, we restrict our attention to the symmetric case. We identify the condition for linear convergence of gradient descent that depends only on the solution $\bm M$ and the sampling set $\Omega$.} In addition, we provide a closed-form expression for the asymptotic convergence rate that matches well with the convergence of the algorithm in practice. \hlnew{The proposed analysis does not require an asymptotic setting for matrix completion, e.g., large matrices of small rank}. We believe that our analysis can be useful in both theoretical and practical aspects of \hlnew{the} matrix completion problem.

\section{Gradient Descent for Matrix Completion}
\label{sec:prel}

\textit{\textbf{Notations.}} Throughout the paper, we use the notations $\norm{\cdot}_F$ and $\norm{\cdot}_2$ to denote the Frobenius norm and the spectral norm of a matrix, respectively. On the other hand, $\norm{\cdot}_2$ is used on a vector to denote the Euclidean norm. Boldfaced symbols are reserved for vectors and matrices. In addition, the $t \times t$ identity matrix is denoted by $\bm I_t$. $\otimes$ denotes the Kronecker product between two matrices, and $\vect(\cdot)$ denotes the vectorization of a matrix  by stacking its columns on top of one another. Let $\bm X$ be some matrix and $\bm F(\bm X)$ be a matrix-valued function of $\bm X$. Then, for some positive number $k$, we use $\bm F(\bm X) = \bm{\mathcal{O}}(\norm{\bm X}_F^k)$ to imply $\lim_{\delta \to 0} \sup_{\norm{\bm X}_F=\delta} \norm{\bm F(\bm X)}_F/\norm{\bm X}_F^k < \infty$.

We begin by introducing the low-rank matrix completion problem. For simplicity, we focus on the {\em symmetric} case where $\bm M$ is an $n \times n$ positive semi-definite (PSD) rank-$r$ matrix and the sampling set $\Omega$ is symmetric.\footnote{If the sampling set is not symmetric, one can symmetrize it by adding $(j,i)$, for any $(i,j) \in \Omega$, to $\Omega$ since $M_{ji} = M_{ij}$.} Let the rank-$r$ economy version of the eigendecomposition of $\bm M$ be given by 
\begin{align*}
    \bm M = \bm U \bm \Lambda \bm U^\topnew,
\end{align*}
where $\bm U \in \R^{n \times r}$ is a semi-orthogonal matrix and $\bm \Lambda \in \R^{r \times r}$ is a diagonal matrix containing $r$ non-zero eigenvalues of $\bm M$, i.e., $\lambda_1 \geq \ldots \geq \lambda_r > 0$.
Since $\bm M$ can be represented as
\begin{align*}
    \bm M = (\bm U \bm \Lambda^{1/2}) (\bm U \bm \Lambda^{1/2})^\topnew ,
\end{align*}
we can write $\bm M = \bm X^* {\bm X^*}^\topnew$, such that $\bm X^* = \bm U \bm \Lambda^{1/2} \in \R^{n \times r}$.
Therefore, the factorization-based \hlnew{formulation} for matrix completion can be \hlnew{described} using the following non-convex optimization:
\begin{align} \label{prob:P0}
    \min_{\bm X \in \R^{n \times r}} \frac{1}{4} \sum_{(i,j) \in \Omega} \bigl( [\bm X \bm X^\topnew]_{ij} - M_{ij} \bigr)^2 .
\end{align}
Denote $\P_\Omega : \R^{n \times n} \to \R^{n \times n}$ the projection onto the set of matrices supported in $\Omega$, i.e.,
\begin{align*}
    [\P_\Omega (\bm Z)]_{ij} = \begin{cases} Z_{ij} & \text{if } (i,j) \in \Omega, \\ 0 & \text{otherwise}. \end{cases}
\end{align*}
We can rewrite the objective function as $f(\bm X) = \frac{1}{4} \norm{\P_\Omega (\bm X \bm X^\topnew - \bm M)}_F^2$.
The gradient of $f(\bm X)$ is given by
\begin{align} \label{equ:df}
    \nabla f(\bm X) = \P_\Omega (\bm X \bm X^\topnew - \bm M) \bm X .
\end{align}
Starting from an initial $\bm X^0$ (usually through spectral initialization \cite{ma2018implicit}), the \hlnew{gradient descent} algorithm (see Algorithm~\ref{algo:GD}) simply updates the value of $\bm X$ by taking steps proportional to the negative of the gradient $\nabla f(\bm X)$.

\begin{algorithm}[t]
\caption{(Non-convex) \hlnew{Gradient Descent}}
\label{algo:GD}
\begin{algorithmic}[1]
\Require{$\bm X^0$, $\P_\Omega(\bm M)$, $\eta$}
\Ensure{$\{\bm X^k\}$}
\For{$k=0,1,2,\ldots$}
\State $\bm X^{k+1} = \bm X^k - \eta \P_\Omega \bigl(\bm X^k {\bm X^k}^\topnew - \bm M\bigr) \bm X^k$
\EndFor
\end{algorithmic}
\end{algorithm}

\section{Local Convergence Analysis}
\label{sec:result}

This section presents the local convergence result of Algorithm~\ref{algo:GD}. While recent work on the global guarantees of the algorithm has shown the linear behavior under certain statistical \hlnew{models}, we emphasize that no closed-form expression of the convergence rate was provided. Our result in this paper, on the other hand, does not make any assumption about the underlying model for $\bm M$ and $\Omega$, and \hlnew{provides} an exact expression of the asymptotic rate of convergence. Let us first introduce some critical concepts used in our derivation.

\begin{definition} \label{def:S}
Denote $\bar{\Omega} = \{ (i-1)n+j \mid (i,j) \in \Omega \}$. 
The row selection matrix $\bm S$ is an ${s \times n^2}$ matrix obtained from a subset of rows corresponding to the elements in $\bar{\Omega}$ from the $n^2 \times n^2$ identity matrix $\bm I_{n^2}$.
\end{definition}

\begin{definition} \label{def:PU}
The orthogonal projection onto the null space of $\bm M$ is defined by $\bm P_{\bm U_\perp} = \bm I_n - \bm U \bm U^\topnew$. 
\end{definition}

\begin{definition} \label{def:T}
Let $\bm T_{n^2}$ be an $n^2 \times n^2$ matrix where the $(i,j)$th block of $\bm T_{n^2}$ is the $n\times n$ matrix $\bm e_j \bm e_i^\topnew$ for $1 \le i,j \le n$. Then $\bm T_{n^2}$ can be used to represent the transpose operator as follows:
\begin{align*}
    \vect(\bm E^\topnew) = \bm T_{n^2} \vect(\bm E) \quad \text{ for any } \bm E \in \R^{n \times n} .
\end{align*}
\end{definition}
We are now in position to state our main result on the asymptotic linear convergence rate of Algorithm~\ref{algo:GD}.
\begin{theorem}
\label{theo:GD}
Denote $\bm P_1 = \bm I_{n^2} - \bm P_{\bm U_\perp} \otimes \bm P_{\bm U_\perp}$, $\bm P_2 = \frac{1}{2} \bigl(\bm I_{n^2} + \bm T_{n^2}\bigr)$, and $\bm P = \bm P_1 \bm P_2$. In addition, let
\begin{align*}
    \bm H = \bm P \Bigl( \bm I_{n^2} - \eta (\bm M \oplus \bm M) (\bm S^\topnew \bm S) \Bigr) \bm P , \numberthis \label{equ:H}
\end{align*}
where $\bm M \oplus \bm M = \bm M \otimes \bm I_n + \bm I_n \otimes \bm M$ is the Kronecker sum. Define the spectral radius of ${\bm H}$, $\rho(\bm H)$, as the largest absolute value of the eigenvalues of $\bm H$.
If $\rho(\bm H)<1$, then Algorithm~\ref{algo:GD} produces a sequence of matrices $\bm X^k {\bm X^k}^\topnew$ converging to $\bm M$ at an asymptotic linear rate $\rho(\bm H)$.
Formally, there exists a neighborhood $\N(\bm M)$ of $\bm M$ such that for any $\bm X^0 {\bm X^0}^\topnew \in \N(\bm M)$,
\begin{align*}
    \norm{\bm X^k {\bm X^k}^\topnew - \bm M}_F \leq C \norm{\bm X^0 {\bm X^0}^\topnew - \bm M}_F \rho(\bm H)^k, \numberthis \label{equ:XXM}
\end{align*}
for some numerical constant $C>0$.
\end{theorem}

\begin{remark}
Theorem~\ref{theo:GD} provides a closed-form expression of the asymptotic linear convergence rate of Algorithm~\ref{algo:GD}, which only depends on $\bm M$, $\Omega$ and the choice of step-size $\eta$. 
We note that the condition for linear convergence, $\rho(\bm H)<1$, is fully determined given $\bm M$, $\Omega$, and $\eta$. 
It would be interesting to establish a connection between this condition and the standard statistical model for matrix completion. For instance, how the incoherence of $\bm M$ and the randomness of $\Omega$ would affect the spectral radius of $\bm H$? 
This exploration is left as future work.
\end{remark}

In our approach, the following lemma plays a pivotal role in the derivation of Theorem~\ref{theo:GD}, establishing the recursion on the error matrix $\bm X^k {\bm X^k}^\topnew - \bm M$:\footnote{We provide proofs of all the lemmas in the Appendix.}

\begin{lemma} \label{lem:EA}
Let $\bm E^k = \bm X^k {\bm X^k}^\topnew - \bm M$. Then
\begin{align*}
    \bm E^{k+1} = \bm E^k - \eta \bigl( \P_\Omega (\bm E^k) \bm M + \bm M \P_\Omega (\bm E^k) \bigr) + \bm{\mathcal{O}}(\norm{\bm E^k}_F^2) .
\end{align*}
Furthermore, denote $\bm A = \bm I_{n^2} - \eta (\bm M \oplus \bm M) (\bm S^\topnew \bm S)$ and $\bm e^k = \vect(\bm E^k)$, the matrix recursion can be rewritten compactly as
\begin{align*}
    \bm e^{k+1} & = \bm A \bm e^k + \bm{\mathcal{O}}(\norm{\bm e^k}_2^2) . \numberthis \label{equ:E4}
\end{align*}
\end{lemma}

\begin{figure}[t]
  \centering
  \centerline{\includegraphics[width=9.1cm]{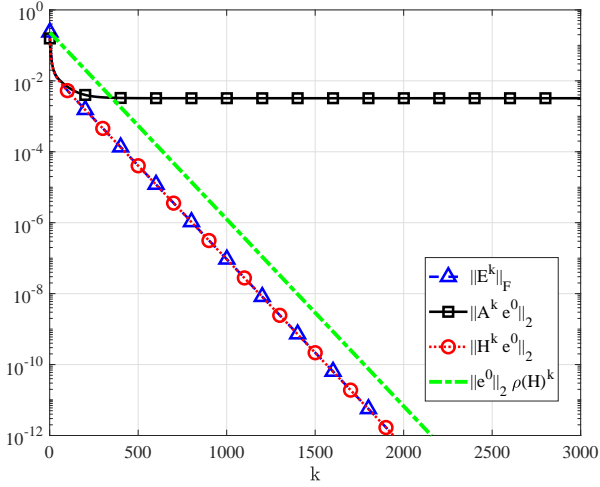}}
  \caption{Linear convergence of gradient descent for matrix completion. \hlnew{The decrease in the norm of error matrix $\bm E^k$ through iterations is shown in the blue dashed line with triangle markers. The black solid line with square markers and the red dotted line with circle markers represent first-order approximations of the error using $\bm A$ and $\bm H$, respectively. Finally, the green dash-dot line is the theoretical bound (up to a constant) given by $\norm{\bm e^0}_2 \rho(\bm H)^k$. We use different markers, i.e., triangle versus circle, to better distinguish the blue line from the red line, respectively.}}
	\label{fig1} 
\vspace{-9pt}
\end{figure}

\begin{remark}
Figure~\ref{fig1} illustrates the effectiveness of the proposed bound on the asymptotic rate of convergence given by Theorem~\ref{theo:GD}. 
\hlnew{In} Fig.~\ref{fig1}, \hlnew{the low-rank solution matrix $\bm M$ is generated by taking the product of a $20 \times 3$ matrix $\bm X$ and its transpose, where $\bm X$ has $i.i.d.$ normally distributed entries. The sampling set $\Omega$ is obtained by randomly selecting the entries of $\bm M$ based on a Bernoulli model with probability $0.3$.
Next, we run the economy-SVD on $\bm M$ to compute $\bm X^* = \bm U \bm \Lambda^{1/2}$. The initialization $\bm X^0$ is obtained by adding $i.i.d.$ normally distributed noise with standard deviation $\sigma=10^{-2}$ to the entries of $\bm X^*$.}
Then we run Algorithm~\ref{algo:GD} \hlnew{with $\bm X^0$, $\P_\Omega(\bm M)$, and $\eta = 0.5/\norm{\bm M}_2$.}
It is noticeable from Fig.~\ref{fig1} that our theoretical bound $\norm{\bm e^0}_2 \rho(\bm H)^k$ \hlnew{given by the green line predicts successfully the rate of decrease in $\norm{\bm E^k}_F$, running parallel to the blue line as soon as $\norm{\bm E^k}_F < 10^{-2}$. As far as the approximations are concerned, we compare the changes in the error modeled by $\bm e^{k+1}=\bm A \bm e^k$ and the error modeled by $\bm e^{k+1}=\bm H \bm e^k$.
While the former (represented by $\norm{\bm A^k \bm e^0}_2$ in black) fails to approximate $\norm{\bm E^k}_F$ for $\norm{\bm E^k}_F < 10^{-2}$, the later (represented by $\norm{\bm H^k \bm e^0}_2$ in red) matches $\norm{\bm E^k}_F$ surprisingly well.} 
\end{remark}

In the rest of this section, we shall derive the proof of Theorem~\ref{theo:GD}. First, we present a major challenge met by the traditional approach that uses (\ref{equ:E4}) to characterize the convergence of the error towards zero. Next, we describe our proposed technique to overcome \hlnew{this difficulty. Finally, we show that our bounding technique recovers the exact rate of local convergence of Algorithm}~\ref{algo:GD}.

\vspace{-5pt}
\subsection{A challenge of establishing the error contraction}
\label{subsec:challenge}
The stability of the nonlinear difference equation (\ref{equ:E4}) is the key to analyze the convergence of Algorithm~\ref{algo:GD}. In essence, linear convergence rate is obtained by the following lemma:

\begin{lemma} 
\label{lem:scalar}
(Rephrased from the supplemental material of \cite{vu2018adaptive}) Let \hlnew{$(a_n)_{n \in \mathbb{N}} \subset \R_+$} be the sequence defined by 
\begin{align*}
a_{n+1} \leq \rho a_{n} + q a_{n}^2 \qquad \text{ for } n=0,1,2,\ldots ,
\end{align*}
where $0\le \rho<1$ and $q \ge 0$. Then $(a_n)$ converges to $0$ if and only if $a_0 < \frac{1-\rho}{q}$. 
A simple linear convergence bound can be derived for $a_0 < \rho \frac{1-\rho}{q}$ in the form of 
\begin{align*}
    a_n \le a_0 K \rho^n, \quad \text{for } K=\biggl(1 - \frac{a_0 q}{ \rho(1-\rho)}\biggr)^{-1} .
\end{align*}
\end{lemma}

\noindent In order to apply Lemma~\ref{lem:scalar} to (\ref{equ:E4}), one natural way is to perform the eigendecomposition $\bm A = \bm Q_{\bm A} \bm \Lambda_{\bm A} \bm Q_{\bm A}^{-1}$, where $\bm Q_{\bm A}$ is the square matrix whose columns are $n^2$ eigenvectors of $\bm A$, and $\bm \Lambda_{\bm A}$ is the diagonal matrix whose diagonal elements are the corresponding eigenvalues of $\bm A$. Then, left-multiplying both sides of (\ref{equ:E4}) by $\bm Q_{\bm A}^{-1}$ yields
\begin{align*}
    \bm Q_{\bm A}^{-1} \bm e^{k+1} & = \bm \Lambda_{\bm A} \bm Q_{\bm A}^{-1} \bm e^k + \bm{\mathcal{O}}(\norm{\bm e^k}_2^2) ,
\end{align*}
\hlnew{where $\bm Q_{\bm A}^{-1}$ does not affect the $\bm{\mathcal{O}}$ term since its norm is constant.}
Applying the triangle inequality\footnote{Given $a=b+c$, by triangle inequality, we have $\|a\| \le \|b\| + \|c\|$ and $\|a\| \ge \|b\| -\|c\|$ (since $b=a+(-c)$ and hence $\|b\| \le \|a\| + \|-c\| = \|a\|+\|c\|$ or $\|a\| \ge \|b\|-\|c\|$). Consequently, we can write $|\|a\| -\|b\|| \le \|c\|$ and hence $\|a\| = \|b\| + {\cal O}(\|c\|)$.} to the last equation leads to 
\begin{align*}
    \norm{\bm Q_{\bm A}^{-1} \bm e^{k+1}}_2 & = \norm{\bm \Lambda_{\bm A} \bm Q_{\bm A}^{-1} \bm e^k}_2 + \mathcal{O}(\norm{\bm e^k}_2^2) . \numberthis \label{equ:eQA}
\end{align*}
With the definition of the spectral radius of $\bm A$ using the spectral norm of $\bm \Lambda_A$, we have
\begin{align*}
    \rho(\bm A) = \norm{\bm \Lambda_{\bm A}}_2 = \sup \biggl\{ \frac{\norm{\bm \Lambda_{\bm A} \tilde{\bm e}}_2}{\norm{\tilde{\bm e}}_2} : \tilde{\bm e} \in \R^{n^2}, \tilde{\bm e} \neq \bm 0 \biggr\} . \numberthis \label{equ:rhoA}
\end{align*}
Now, using (\ref{equ:rhoA}) and the fact that $\mathcal{O}(\norm{\bm e^k}_2^2) = \mathcal{O}(\norm{\bm Q_{\bm A}^{-1} \bm e^k}_2^2)$, (\ref{equ:eQA}) can be upper-bounded by
\begin{align*}
    \norm{\bm Q_{\bm A}^{-1} \bm e^{k+1}}_2 & \leq \rho(\bm A) \norm{\bm Q_{\bm A}^{-1} \bm e^k}_2 + \mathcal{O}(\norm{\bm Q_{\bm A}^{-1} \bm e^k}_2^2) . \numberthis \label{equ:boundA}
\end{align*}
If $\rho(\bm A) < 1$, then by Lemma~\ref{lem:scalar}, the sequence $\norm{\bm Q_{\bm A}^{-1} \bm e^k}_2$ converges to $0$ linearly at rate $\rho(\bm A)$. Unfortunately, one can verify that $\rho(\bm A) \geq 1$ by taking any vector $\bm v \in \R^{n^2}$ such that $v_i=0$ for all $i \in \bar{\Omega}$. Since $\bm A \bm v = \bm v$, $1$ must be an eigenvalue of $\bm A$. 

The failure of the aforementioned bounding technique is it overlooks the fact that $\bm E^k = \bm X^k {\bm X^k}^\topnew - \bm M$. By defining $\E = \{ \bm X \bm X^\topnew - \bm M \mid \bm X \in \R^{n \times r} \}$ and $\tilde{\E}_{\bm A} = \{ \bm Q^{-1}_{\bm A} \vect(\bm E) \mid \bm E \in \E \}$,
a tighter bound \hlnew{on $\norm{\bm \Lambda_{\bm A} \bm Q_{\bm A}^{-1} \bm e^k}_2 / \norm{\bm Q_{\bm A}^{-1} \bm e^k}_2$} can be obtained by
\begin{align} \label{equ:rho_constrained}
    \rho^{\E}(\bm A,\delta) = \sup \biggl\{ \frac{\norm{\bm \Lambda_{\bm A} \tilde{\bm e}}_2}{\norm{\tilde{\bm e}}_2} : \tilde{\bm e} \in \tilde{\E}_{\bm A}, \tilde{\bm e} \neq \bm 0, \norm{\tilde{\bm e}}_2 \leq \delta \biggr\} ,
\end{align}
for some constant $\delta > 0$.
Taking into account the structure of $\bm E^k$, one would expect $\rho^\E (\bm A) = \lim_{\delta \to 0} \rho^{\E}(\bm A,\delta)$ is a more reliable estimate of the asymptotic rate of convergence for (\ref{equ:E4}). Nonetheless, (\ref{equ:rho_constrained}) is a non-trivial optimization problem that has no closed-form solution to the best of our knowledge.

\subsection{Integrating structural constraints}

\hlnew{To address the aforementioned issue, we propose to integrate the structural constraint on $\bm E^k$ into the recursion} (\ref{equ:E4}). \hlnew{As we shall show in the next subsection, this integration enables the application of Lemma}~\ref{lem:scalar} \hlnew{to the new recursion in order to obtain a tight bound on the convergence rate. First, let us characterize the feasible set of error matrices $\E$ as follows:}
\begin{lemma} \label{lem:E}
$\bm E \in \E$ if and only if the following conditions hold simultaneously:
\begin{enumerate}
    \item[(C1)] $\P_r(\bm M + \bm E) = \bm M + \bm E$, where $\P_r$ is the truncated singular value decomposition of order $r$ \cite{eckart1936approximation}.
    \item[(C2)] $\bm E^\topnew = \bm E$.
    \item[(C3)] $\bm v^\topnew (\bm M + \bm E) \bm v \geq 0$ for all $\bm v \in \R^n$.
\end{enumerate}
\end{lemma}
\noindent 
Our strategy is to integrate three conditions in Lemma~\ref{lem:E} into the linear operator $\bm A$ so that the resulting recursion will implicitly enforce $\bm E^k$ to remain in $\E$. Specifically, for condition (C1), we linearize $\P_r$ using the first-order perturbation analysis of the truncated singular value decomposition \cite{vu2020perturbation}. For condition (C2), we leverage the linearity of the transpose operator. \hlnew{Finally, while handling condition (C3) is non-trivial, it turns out that this condition can be ignored.}
In the following lemma, we introduce the linear projection that ensures the updated error ${\bm E}^k$ remains near $\E$.
\begin{lemma} \label{lem:e}
Recall that $\bm P_1 = \bm I_{n^2} - \bm P_{\bm U_\perp} \otimes \bm P_{\bm U_\perp}$, $\bm P_2 = \frac{1}{2} \bigl(\bm I_{n^2} + \bm T_{n^2}\bigr)$. Then, the following statements hold:
\begin{enumerate}
    \item $\bm P_1$ corresponds to the orthogonal projection onto the tangent plane of the set of rank-$r$ matrices at $\bm M$.
    \item $\bm P_2$ corresponds to the orthogonal projection onto the space of symmetric matrices.
    \item $\bm P_1$ and $\bm P_2$ commute, and $\bm P = \bm P_1 \bm P_2 = \bm P_2 \bm P_1$ is also an orthogonal projection.
    \item For any $\bm E \in \E$, $\vect (\bm E) = \bm P \vect (\bm E) + \bm{\mathcal{O}}(\norm{\bm E}_F^2)$.
\end{enumerate}
\end{lemma}
\noindent By Lemma~\ref{lem:e}-4, we have $\bm e^{k} = \bm P \bm e^{k} + \bm{\mathcal{O}}(\norm{\bm e^{k}}_2^2)$ for all $k$. Using this result with $k+1$ instead of $k$ and replacing $\bm e^{k+1}$ from (\ref{equ:E4}) into the first term on the RHS, we have $$\bm e^{k+1} = \bm P \bigl(\bm A \bm e^k + \bm{\mathcal{O}}(\norm{\bm e^k}_2^2)\bigr) + \bm{\mathcal{O}}(\norm{\bm e^{k+1}}_2^2) .$$
Substituting $\bm e^{k} = \bm P \bm e^{k} + \bm{\mathcal{O}}(\norm{\bm e^{k}}_2^2)$ and using $\bm e^{k+1} = \bm{\mathcal{O}}(\norm{\bm e^k}_2)$, we obtain
\begin{align*}
    \bm e^{k+1} &= \bm P \bm A \bm P \bm e^k + \bm{\mathcal{O}}(\norm{\bm e^k}_2^2) . \numberthis \label{equ:E6}
\end{align*}
It can be seen from Lemma~\ref{lem:e}-1 and Lemma~\ref{lem:e}-2 that the projection $\bm P$ enforces the error vector $\bm e^k$ to lie in the space under conditions (C1) and (C2) in Lemma~\ref{lem:E}.
Now replacing the definition $\bm H = \bm P \bm A \bm P$, (\ref{equ:E6}) can be rewritten as
\begin{align}
    \bm e^{k+1} = \bm H \bm e^k + \bm{\mathcal{O}}(\norm{\bm e^k}_2^2) . \numberthis \label{equ:E7}
\end{align}
Similar to the derivation with $\bm A$, let $\bm H = \bm Q_{\bm H} \bm \Lambda_{\bm H} \bm Q_{\bm H}^{-1}$ be the eigendecomposition of $\bm H$ and define $\tilde{\bm e}^k = \bm Q^{-1}_{\bm H} \bm e^k$. Then, we have
\begin{align*}
    \norm{\tilde{\bm e}^{k+1}}_2 = \norm{\bm \Lambda_{\bm H} \tilde{\bm e}^k}_2 + \mathcal{O}(\norm{\tilde{\bm e}^k}_2^2) . \numberthis \label{equ:E8}
\end{align*}
In addition, denote $\tilde{\E}_{\bm H} = \{ \bm Q^{-1}_{\bm H} \vect(\bm E) \mid \bm E \in \E \}$, we can define
\begin{align*}
    \rho(\bm H) &= \sup \biggl\{ \frac{\norm{\bm \Lambda_{\bm H} \tilde{\bm e}}_2}{\norm{\tilde{\bm e}}_2} : \tilde{\bm e} \in \R^{n^2}, \tilde{\bm e} \neq \bm 0 \biggr\} \text{ and } \numberthis \label{equ:rhoH} \\ 
    \rho^{\E}(\bm H,\delta) &= \sup \biggl\{ \frac{\norm{\bm \Lambda_{\bm H} \tilde{\bm e}}_2}{\norm{\tilde{\bm e}}_2} : \tilde{\bm e} \in \tilde{\E}_{\bm H}, \tilde{\bm e} \neq \bm 0, \norm{\tilde{\bm e}}_2 \leq \delta \biggr\} . \numberthis \label{equ:rho_constrainedH}
\end{align*}
Since (\ref{equ:E4}) and (\ref{equ:E7}) are two different systems that describes the same dynamic for $\bm E^k \in \E$, one would expect they share the same asymptotic behavior. In particular, their linear rates of convergence should agree when the constraint $\bm E^k \in \E$ is considered.
\begin{lemma} \label{lem:rho}
Let $\rho^\E (\bm H) = \lim_{\delta \to 0} \rho^{\E}(\bm H,\delta)$. Then, 
\begin{align*}
    \rho^\E (\bm H) = \rho^\E (\bm A).
\end{align*}
\end{lemma}
\noindent \hlnew{While using $\bm H$ instead of $\bm A$ preserves the system dynamic over $\E$, it provides updates of the error that ensure that it remains in $\E$. Consequently, we can ignore the constraints that are implicitly satisfied in our analysis when using $\bm H$.}

\subsection{Asymptotic bound on the linear convergence rate}
\hlnew{We have seen in Subsection}~\ref{subsec:challenge} that applying Lemma~\ref{lem:scalar} to (\ref{equ:E4}) \hlnew{fails to estimate the convergence rate due to the gap between $\rho^\E(\bm A)$ and $\rho(\bm A)$. In this subsection, we show that integrating the structural constraint helps eliminating the gap between $\rho^\E(\bm H)$ and $\rho(\bm H)$ (even when condition (C3) is omitted). Therefore, applying} Lemma~\ref{lem:scalar} to (\ref{equ:E8}) yields $\rho(\bm H)$ as a tight bound on the convergence rate. To that end, our goal is to prove the following lemma:
 \begin{lemma} \label{lem:tight}
As $\delta$ approaches $0$, we have $\rho(\bm H) - \rho^{\E}(\bm H,\delta) = \mathcal{O}(\delta)$. Consequently, it holds that $\rho(\bm H) = \rho^{\E}(\bm H)$.
\end{lemma}
\noindent \hlnew{Let us briefly present the key ideas and lemmas we use to prove Lemma}~\ref{lem:tight}.
\hlnew{Our proof relies on two critical considerations: (i) $\rho^\E (\bm H,\delta) \leq \rho(\bm H)$, (ii) there exists a maximizer $\tilde{\bm e}^\star$ of the supremum in} (\ref{equ:rhoH}) such that the distance from $\tilde{\bm e}^\star$ to $\tilde{\E}_{\bm H}$ is $\mathcal{O}(\delta^2)$. While (i) is trivial from (\ref{equ:rhoH}) and (\ref{equ:rho_constrainedH}), \hlnew{(ii) is proven by introducing $\F_\delta$ as a surrogate for the set $\E$ as follows:}
\begin{lemma} \label{lem:PSD}
\hlnew{Denote the eigenvector of $\bm H$ corresponding to the largest (in magnitude) eigenvalue by $\bm q_1$. Define $\bm G$ as the ${n \times n}$ matrix satisfying $\vect(\bm G) = \delta \bm q_1$.}
Let $\F_\delta$ be the set of $n \times n$ matrices satisfying the following conditions: (i) $\norm{\bm F}_F \leq 2\delta$; (ii) $\bm F^\topnew = \bm F$; (iii) $\norm{\bm P_{\bm U_\perp} \bm F \bm P_{\bm U_\perp}}_F \leq \frac{2}{\lambda_r} \delta^2$; and (iv) $\bm v^\topnew (\bm M + \bm F) \bm v \geq 0$ for all $\bm v \in \R^n$.
Then, there exists $\bm F \in \F_\delta$ satisfying
\begin{align*}
    \norm{\bm F - \bm G}_F = \mathcal{O}(\delta^2) .
\end{align*} 
\end{lemma}
\begin{lemma} \label{lem:Pr}
For any $\bm F \in \F_\delta$, there exists $\bm E \in \E$ satisfying
\begin{align*}
    \norm{\bm E - \bm F}_F = \mathcal{O}(\delta^2) .
\end{align*}
\end{lemma}
\noindent \hlnew{From (i) and (ii), it follows that the difference between $\rho^\E(\bm H,\delta)$ and $\rho(\bm H)$ is $\mathcal{O}(\delta)$. Thus, $\rho(\bm H) = \rho^\E(\bm H)$ when taking the limit of $\rho^\E(\bm H,\delta)$ as $\delta \to 0$. Our derivation of Theorem}~\ref{theo:GD} is completed by directly applying Lemma~\ref{lem:scalar} to (\ref{equ:E8}). 

\section{Conclusion and Future work}
\label{sec:conc}
We presented a framework for analyzing the convergence of the existing gradient descent approach for low-rank matrix completion. In our analysis, we restricted our focus to the symmetric matrix completion case. We proved that the algorithm converges linearly. Different to other approaches, we made no assumption on the rank of the matrix or fraction of available entries. Instead, we derived an expression for the linear convergence rate via the spectral norm of a closed-form matrix. As future work, using random matrix theory, the closed-form expression for the convergence rate can be further related to the rank, the number of available entries, and the matrix dimensions. Additionally, this work can be extended to the non-symmetric case.

\section{Appendix}

\subsection{Proof of Lemma~\ref{lem:EA}}
Recall the gradient descent update in Algorithm~\ref{algo:GD}:
\begin{align*} 
    \bm X^{k+1} &= \bm X^k - \eta \P_\Omega \bigl(\bm X^k {\bm X^k}^\topnew - \bm M\bigr) \bm X^k \\
    &= (\bm I_n  - \eta \P_\Omega (\bm E^k)) \bm X^k . \numberthis \label{equ:GD}
\end{align*}
Substituting (\ref{equ:GD}) into the definition of $\bm E^{k+1}$, we have
\begin{align*}
    &\bm E^{k+1} = \bm X^{k+1} {\bm X^{k+1}}^\topnew - \bm M \\
    &\quad = \bigl( \bm I_n  - \eta \P_\Omega (\bm E^k)  \bigr) \bm X^k {\bm X^k}^\topnew \bigl( \bm I_n  - \eta \P_\Omega (\bm E^k)  \bigr)^\topnew - \bm M.
\end{align*}
From the fact that $\bm E^k$ is symmetric and $\Omega$ is a symmetric sampling, the last equation can be further expanded as
\begin{align*}
    &\bm E^{k+1} = \bm X^k {\bm X^k}^\topnew - \eta \P_\Omega (\bm E^k) \bm X^k {\bm X^k}^\topnew \\
    &- \eta \bm X^k {\bm X^k}^\topnew \P_\Omega (\bm E^k) + \eta^2 \P_\Omega (\bm E^k) \bm X^k {\bm X^k}^\topnew \P_\Omega (\bm E^k) - \bm M . \numberthis \label{equ:E1}
\end{align*}
Since $\bm X^k {\bm X^k}^\topnew = \bm M + \bm E^k$, (\ref{equ:E1}) is equivalent to
\begin{align*}
    \bm E^{k+1} &= \bm E^k - \eta \bigl( \P_\Omega (\bm E^k) \bm M + \bm M \P_\Omega (\bm E^k) \bigr) \\
    &\quad - \eta \bigl( \P_\Omega (\bm E^k) \bm E^k + \bm E^k \P_\Omega (\bm E^k) \bigr) \\
    &\quad + \eta^2 \P_\Omega (\bm E^k) \bm M \P_\Omega (\bm E^k) + \eta^2 \P_\Omega (\bm E^k) \bm E^k \P_\Omega (\bm E^k) . \numberthis \label{equ:E2}
\end{align*}
Note that $\norm{\P_\Omega (\bm E^k)}_F \leq \norm{\bm E^k}_F$. Hence, collecting terms that are of second order and higher, with respect to $\norm{\bm E^k}_F$, on the RHS of (\ref{equ:E2}) yields
\begin{align*}
    \bm E^{k+1} &= \bm E^k - \eta \bigl( \P_\Omega (\bm E^k) \bm M + \bm M \P_\Omega (\bm E^k) \bigr) + \bm{\mathcal{O}}(\norm{\bm E^k}_F^2) .
\end{align*}
Now by Definition~\ref{def:S}, it is easy to verify that
$$ \bm S \bm S^\topnew = \bm I_{n^2} \quad \text{and} \quad \vect\bigl(\P_\Omega (\bm E^k)\bigr) = \bm S^\topnew \bm S \bm e^k .$$
Using the property 
$\vect(\bm A \bm B \bm C) = (\bm C^\topnew \otimes \bm A) \vect(\bm B)$, (5) can be vectorized as follows: 
\begin{align*}
    \bm e^{k+1} &= \bm e^k - \eta (\bm M \otimes \bm I_n) \vect\bigl(\P_\Omega (\bm E^k)\bigr) \\
    &\quad - \eta (\bm I_n \otimes \bm M) \vect\bigl(\P_\Omega (\bm E^k)\bigr) + \bm{\mathcal{O}}(\norm{\bm e^k}_2^2) .
\end{align*}
The last equation can be reorganized as
\begin{align*}
    \bm e^{k+1} = \Bigl( \bm I_{n^2} - \eta (\bm M \oplus \bm M) (\bm S^\topnew \bm S) \Bigr) \bm e^k + \bm{\mathcal{O}}(\norm{\bm e^k}_2^2) .
\end{align*}


\subsection{Proof of Lemma~\ref{lem:E}}
($\Rightarrow$) Suppose $\bm E \in \E$. Then for (C1), i.e., $\bm E^\topnew = \bm E$, $\bm E = \bm X \bm X^\topnew - \bm M$ is symmetric since both $\bm X \bm X^\topnew$ and $\bm M$ are symmetric. 
For (C2), i.e., $\P_r (\bm M + \bm E) = \bm M + \bm E$, stems from the fact $\bm M + \bm E = \bm X \bm X^\topnew$ has rank no greater than $r$ for $\bm X \in \R^{n \times r}$. Finally, for any $\bm v \in \R^n$, we have
\begin{align*}
    \bm v^\topnew (\bm M + \bm E) \bm v = \bm v^\topnew (\bm X \bm X^\topnew) \bm v = \norm{\bm X^\topnew \bm v}_2^2 \geq 0 .
\end{align*}

\noindent ($\Leftarrow$) From conditions (C1) and (C3), $\bm M + \bm E$ is a PSD matrix. In addition, $\P_r(\bm M + \bm E) = \bm M + \bm E$ implies $\bm M + \bm E$ must have rank no greater $r$. Since any PSD matrix $\bm A$ with rank less than or equal to $r$ can be factorized as $\bm A = \bm Y \bm Y^\topnew$ for some $\bm Y \in \R^{n \times r}$, we conclude that $\bm E \in \E$.

\subsection{Proof of Lemma~\ref{lem:e}}
First, recall that any matrix $\bm \Pi \in \R^{n^2 \times n^2}$ is an orthogonal projection if and only if $\bm \Pi^2 = \bm \Pi$ and $\bm \Pi = \bm \Pi^\topnew$. Since $\bm P_{\bm U_\perp}^\topnew = \bm P_{\bm U_\perp}$, we have
\begin{align*}
    \bm P_1^\topnew &= \bigl(\bm I_{n^2} - \bm P_{\bm U_\perp} \otimes \bm P_{\bm U_\perp}\bigr)^\topnew \\
    &= \bm I_{n^2}^\topnew - \bm P_{\bm U_\perp}^\topnew \otimes \bm P_{\bm U_\perp}^\topnew \\
    &= \bm I_{n^2} - \bm P_{\bm U_\perp} \otimes \bm P_{\bm U_\perp} = \bm P_1 .
\end{align*}
In addition, since $\bm P_{\bm U_\perp}^2 = \bm P_{\bm U_\perp}$, we have
\begin{align*}
    \bm P_1 ^2 &= (\bm I_{n^2} - \bm P_{\bm U_\perp} \otimes \bm P_{\bm U_\perp}) (\bm I_{n^2} - \bm P_{\bm U_\perp} \otimes \bm P_{\bm U_\perp})^\topnew \\
    &= \bm I_{n^2}^2 - 2 \bm P_{\bm U_\perp} \otimes \bm P_{\bm U_\perp} + (\bm P_{\bm U_\perp} \otimes \bm P_{\bm U_\perp})^2 \\
    &= \bm I_{n^2} - 2 \bm P_{\bm U_\perp} \otimes \bm P_{\bm U_\perp} + (\bm P_{\bm U_\perp}^2 \otimes \bm P_{\bm U_\perp}^2) \\
    &= \bm I_{n^2} - 2 \bm P_{\bm U_\perp} \otimes \bm P_{\bm U_\perp} + \bm P_{\bm U_\perp} \otimes \bm P_{\bm U_\perp} \\
    &= \bm I_{n^2} - \bm P_{\bm U_\perp} \otimes \bm P_{\bm U_\perp} = \bm P_1.
\end{align*}
Second, using the fact that $\bm T_{n^2}^2 = \bm I_{n^2}$ and $\bm T_{n^2}$ is symmetric, we can derive similar result: 
\begin{align*}
    \bm P_2^\topnew = \biggl( \frac{\bm I_{n^2} + \bm T_{n^2}}{2} \biggr)^\topnew = \frac{\bm I_{n^2} + \bm T_{n^2}}{2} = \bm P_2 ,
\end{align*}
and 
\begin{align*}
    \bm P_2^2 &= \frac{(\bm I_{n^2} + \bm T_{n^2})^2}{4} \\
    &= \frac{\bm I_{n^2} + 2 \bm T_{n^2} + \bm T^2_{n^2}}{4} \\
    &= \frac{2\bm I_{n^2} + 2\bm T_{n^2}}{4} \\
    &= \frac{\bm I_{n^2} + \bm T_{n^2}}{2} = \bm P_2 .
\end{align*}

Third, we observe that $\bm P_1$ and $\bm P_2$ are the vectorized version of the linear operators
\begin{align*}
    \bm \Pi_1 (\bm E) = \bm E - \bm P_{\bm U_\perp} \bm E \bm P_{\bm U_\perp}
\end{align*}
and
\begin{align*}
    \bm \Pi_2 (\bm E) = \frac{1}{2} (\bm E + \bm E^\topnew) ,
\end{align*}
respectively, for any $\bm E \in \R^{n \times n}$. Hence, in order to prove that $\bm P_1$ and $\bm P_2$ commute, it is sufficient to show that operators $\bm \Pi_1$ and $\bm \Pi_2$ commute. Indeed, we have
\begin{align*}
    \bm \Pi_2 \bm \Pi_1 (\bm E) &= \frac{1}{2} \bigl( (\bm E - \bm P_{\bm U_\perp} \bm E \bm P_{\bm U_\perp}) + (\bm E - \bm P_{\bm U_\perp} \bm E \bm P_{\bm U_\perp})^\topnew \bigr) \\
    &= \frac{1}{2} (\bm E + \bm E^\topnew) - \bm P_{\bm U_\perp} \frac{1}{2} (\bm E + \bm E^\topnew) \bm P_{\bm U_\perp} \\
    &= \bm \Pi_1 \bm \Pi_2 (\bm E) .
\end{align*}
This implies $\bm \Pi_1$ and $\bm \Pi_2$ commute.
Since $\bm P$ is the product of two commuting orthogonal projections, it is also an orthogonal projection.

Finally, let us restrict $\bm E$ to belong to $\E$ and denote $\bm e = \vect (\bm E)$. Using Theorem~3 in \cite{vu2020perturbation}, we have
\begin{align*}
    \P_r (\bm M + \bm E) = \bm M + \bm E - \bm P_{\bm U_\perp} \bm E \bm P_{\bm U_\perp} + \bm{\mathcal{O}}(\norm{\bm E}_F^2) . \numberthis \label{equ:Pr}
\end{align*}
Since $\P_r (\bm M + \bm E) = \bm M + \bm E$, it follows from (\ref{equ:Pr}) that 
\begin{align*}
    \bm P_{\bm U_\perp} \bm E \bm P_{\bm U_\perp} = \bm{\mathcal{O}}(\norm{\bm E}_F^2) .
\end{align*}
Vectorizing the last equation, we obtain 
\begin{align} \label{equ:P1}
    (\bm P_{\bm U_\perp} \otimes \bm P_{\bm U_\perp}) \bm e = \bm{\mathcal{O}}(\norm{\bm E}_F^2) .
\end{align}
On the other hand, since $\bm E$ is symmetric,
\begin{align} \label{equ:P2}
    \bm e = \bm T_{n^2} \bm e = \Bigl( \frac{\bm I_{n^2} + \bm T_{n^2}}{2} \Bigr) \bm e .
\end{align}
From (\ref{equ:P1}) and (\ref{equ:P2}), we have
\begin{align*}
    \bm e &= (\bm I_{n^2} - \bm P_{\bm U_\perp} \otimes \bm P_{\bm U_\perp}) \bm e + \bm{\mathcal{O}}(\norm{\bm E}_F^2) \\
    &= (\bm I_{n^2} - \bm P_{\bm U_\perp} \otimes \bm P_{\bm U_\perp}) \Bigl( \frac{\bm I_{n^2} + \bm T_{n^2}}{2} \Bigr) \bm e + \bm{\mathcal{O}}(\norm{\bm E}_F^2) . \numberthis \label{equ:e_recur}
\end{align*}
Substituting 
$$\bm P = \bm P_1 \bm P_2 = (\bm I_{n^2} - \bm P_{\bm U_\perp} \otimes \bm P_{\bm U_\perp}) \Bigl( \frac{\bm I_{n^2} + \bm T_{n^2}}{2} \Bigr) $$
into (\ref{equ:e_recur}) completes our proof of the lemma.

\subsection{Proof of Lemma~\ref{lem:rho}}
Let $\tilde{\E} = \{ \vect(\bm E) \mid \bm E \in \E \}$.
Recall that for any $\bm e \in \tilde{\E}$,
\begin{align*}
    \bm e = \bm P \bm e + \bm{\mathcal{O}}(\norm{\bm e}_2^2) .
\end{align*}
Therefore, by the triangle inequality, we obtain
\begin{align*}
    \norm{\bm A \bm e}_2 &= \norm{\bm A \bigl(\bm P \bm e + \bm{\mathcal{O}}(\norm{\bm e}_2^2) \bigr)} \\
    &\leq \norm{\bm A \bm P \bm e}_2 + \norm{\bm A \bm{\mathcal{O}}(\norm{\bm e}_2^2)}_2 .
\end{align*}
Since the second term on the RHS of the last inequality is $\mathcal{O}(\norm{\bm e}_2^2)$, it is also $\mathcal{O}(\delta^2)$ for any $\bm e \in \tilde{\E}$ such that $\norm{\bm e}_2 \leq \delta$. In other words,
\begin{align} \label{equ:uAe}
    \norm{\bm A \bm e}_2 \leq \norm{\bm A \bm P \bm e}_2 + \mathcal{O}(\delta^2) .
\end{align}
Similarly, we also have,
\begin{align*} 
    \norm{\bm A \bm e}_2 &\geq \norm{\bm A \bm P \bm e}_2 - \norm{\bm A \bm{\mathcal{O}}(\norm{\bm e}_2^2)}_2 \\
    &= \norm{\bm A \bm P \bm e}_2 - \mathcal{O}(\delta^2) . \numberthis \label{equ:lAe}
\end{align*}
From (\ref{equ:uAe}) and (\ref{equ:lAe}), it follows that
\begin{align*} 
    \frac{\norm{\bm A \bm e}_2}{\norm{\bm e}_2} &= \frac{\norm{\bm A \bm P \bm e}_2}{\norm{\bm e}_2} + \mathcal{O}(\delta) . \numberthis \label{equ:ulAe}
\end{align*}
Taking the limit of the supremum of (\ref{equ:ulAe}) as $\delta \to 0$ yields
\begin{align*}
    \rho^\E(\bm A) &= \lim_{\delta \to 0} \sup_{\substack{\bm e \in \tilde{\E} \\ \bm e \neq 0 \\ \norm{\bm e}_2 \leq \delta}} \frac{\norm{\bm A \bm e}_2}{\norm{\bm e}_2} \\
    &= \lim_{\delta \to 0} \sup_{\substack{\bm e \in \tilde{\E} \\ \bm e \neq 0 \\ \norm{\bm e}_2 \leq \delta}} \frac{\norm{\bm A \bm P \bm e}_2}{\norm{\bm e}_2} = \rho^\E(\bm A \bm P) . \numberthis \label{equ:A_AP}
\end{align*}
Now following similar argument in Lemma~\ref{lem:tight}, we have
\begin{align} \label{equ:PAP}
    \begin{cases} \rho^\E(\bm A \bm P) = \rho(\bm A \bm P) , \\
    \rho^\E(\bm P \bm A \bm P) = \rho(\bm P \bm A \bm P) .
    \end{cases}
\end{align}
Given (\ref{equ:A_AP}) and (\ref{equ:PAP}), it remains to show that $\rho(\bm A \bm P) = \rho(\bm P \bm A \bm P)$. Indeed, using Gelfand's formula \cite{gelfand1941normierte}, we have
\begin{align*}
    &\rho(\bm A \bm P) = \lim_{k \to \infty} \norm{(\bm A \bm P)^k}_2^{1/k} \\
    \text{ and } &\rho(\bm P\bm A \bm P) = \lim_{k \to \infty} \norm{(\bm P \bm A \bm P)^k}_2^{1/k} .
\end{align*}
By the property of operator norms,
\begin{align*}
    \norm{(\bm A \bm P)^k}_2 = \norm{\bm A (\bm P \bm A \bm P)^{k-1}}_2 \leq \norm{\bm A}_2 \norm{(\bm P \bm A \bm P)^{k-1}}_2 .
\end{align*}
Thus,
\begin{align*}
    \norm{(\bm A \bm P)^k}_2^{1/k} \leq \norm{\bm A}_2^{1/k} \Bigl(\norm{(\bm P \bm A \bm P)^{k-1}}_2^{1/(k-1)} \Bigr)^{(k-1)/k} .
\end{align*}
Taking the limit of both sides of the last inequality as $k \to \infty$ yields $\rho(\bm A \bm P) \leq \rho(\bm P \bm A \bm P)$.
Similarly, since
\begin{align*}
    \norm{(\bm P \bm A \bm P)^k}_2 = \norm{\bm P (\bm A \bm P)^k}_2 \leq \norm{(\bm A \bm P)^k}_2 ,
\end{align*}
we also obtain $\rho(\bm P \bm A \bm P) \leq \rho(\bm A \bm P)$. This concludes our proof of the lemma.

\subsection{Proof of Lemma~\ref{lem:tight}}

Without loss of generality, assume $\lambda_1$ is the eigenvalue with largest magnitude, i.e., $\abs{\lambda_1} = \rho(\bm H)$. By the definition of $\bm G$, we have $\norm{\bm G}_F = \delta$. Since $\bm H \vect(\bm G) = \lambda_1 \vect(\bm G)$ and $\bm H = \bm Q_{\bm H} \bm \Lambda_{\bm H} \bm Q^{-1}_{\bm H}$, it follows that
\begin{align} \label{equ:QHG}
    \bm Q_{\bm H} \bm \Lambda_{\bm H} \bm Q_{\bm H}^{-1} \vect(\bm G) = \lambda_1 \vect(\bm G).
\end{align}
Multiplying both sides of (\ref{equ:QHG}) by $\bm Q_{\bm H}^{-1}$, we obtain
\begin{align*}
    \bm \Lambda_{\bm H} \bm Q_{\bm H}^{-1} \vect(\bm G) = \lambda_1 \bm Q_{\bm H}^{-1}\vect(\bm G) . 
\end{align*}
Taking the $L2$-norm and and reorganizing the equation yields
\begin{align*}
     \frac{\norm{\bm \Lambda_{\bm H} \bm Q_{\bm H}^{-1} \vect(\bm G) }_2}{\norm{\bm Q_{\bm H}^{-1} \vect(\bm G)}_2} = \abs{\lambda_1} = \rho(\bm H) . \numberthis \label{equ:Grho}
\end{align*}
Therefore, $\bm G$ leads to a solution of the supremum in (\ref{equ:rhoH}). We now prove that $\bm G$ is symmetric and $(\bm P_{\bm U_\perp} \otimes \bm P_{\bm U_\perp}) \vect(\bm G) = \bm 0$.
First, since $\bm P_1$, $\bm P_2$ and $\bm P = \bm P_1 \bm P_2$ are orthogonal projections, we have
\begin{align*}
    \bm P_2 \bm H &= \bm P_2 \bm P \bm A \bm P \\
    &= \bm P_2 \bm P_2 \bm P_1 \bm A \bm P \\
    &= \bm P_2 \bm P_1 \bm A \bm P \\
    &= \bm P_1 \bm P_2 \bm A \bm P \\
    &= \bm P \bm A \bm P = \bm H .
\end{align*}
Thus, 
\begin{align*}
    \lambda_1 \vect(\bm G) &= \bm H \vect(\bm G) \\
    &= \bm P_2 \bm H \vect(\bm G) \\
    &= \lambda_1 \bm P_2 \vect(\bm G) . \numberthis \label{equ:Gt}
\end{align*}
Substituting $\bm P_2 = \frac{1}{2} \bigl(\bm I_{n^2} + \bm T_{n^2}\bigr)$ into (\ref{equ:Gt}) yields
\begin{align*}
    \vect(\bm G^\topnew) = \bm T_{n^2} \vect(\bm G) \quad \text{ or } \quad \bm G = \bm G^\topnew .
\end{align*}
Second, since $\bm P_1 \bm H = \bm H$, we obtain
\begin{align*}
    \lambda_1 \vect(\bm G) &= \bm H \vect(\bm G) \\
    &= \bm P_1 \bm H \vect(\bm G) \\
    &= \lambda_1 \bm P_1 \vect(\bm G) . \numberthis \label{equ:GP}
\end{align*}
Substituting $\bm P_1 = \bm I_{n^2} - \bm P_{\bm U_\perp} \otimes \bm P_{\bm U_\perp}$ into (\ref{equ:GP}) yields
\begin{align*}
    (\bm P_{\bm U_\perp} \otimes \bm P_{\bm U_\perp}) \vect(\bm G) = \bm 0 \quad \text{ or } \quad \bm P_{\bm U_\perp} \bm G \bm P_{\bm U_\perp} = \bm 0 .
\end{align*}
Since $\norm{\bm E - \bm G}_F \leq \norm{\bm E - \bm F}_F + \norm{\bm F - \bm G}_F$ (by the triangle inequality), Lemmas~\ref{lem:PSD} and \ref{lem:Pr} imply the existence of $\bm E \in \E$ such that $\norm{\bm E - \bm G}_F = \mathcal{O}(\delta^2)$. 
Denote $\tilde{\bm e} = \bm Q^{-1}_{\bm H} \vect(\bm E) \in \tilde{\E}_{\bm H}$, we have
\begin{align*}
    \bm \Lambda_{\bm H} \tilde{\bm e} &= \lambda_1 \tilde{\bm e} - (\lambda_1 \bm I_{n^2} - \bm \Lambda_{\bm H}) \tilde{\bm e} \\
    &= \lambda_1 \tilde{\bm e} - (\lambda_1 \bm I_{n^2} - \bm \Lambda_{\bm H}) \bm Q^{-1}_{\bm H} \vect(\bm E) \\
    &= \lambda_1 \tilde{\bm e} - (\lambda_1 \bm I_{n^2} - \bm \Lambda_{\bm H}) \bm Q^{-1}_{\bm H} \vect(\bm E - \bm G) .
\end{align*}
where the last equality stems from the fact that $\lambda_1 \bm Q^{-1}_{\bm H} \vect(\bm G) = \bm \Lambda_{\bm H} \bm Q^{-1}_{\bm H} \vect(\bm G)$.
Next, using the triangle inequality, we obtain
\begin{align*}
    &\norm{\bm \Lambda_{\bm H} \tilde{\bm e}}_2 \geq \norm{\lambda_1 \tilde{\bm e}}_2 - \norm{(\lambda_1 \bm I_{n^2} - \bm \Lambda_{\bm H}) \bm Q^{-1}_{\bm H} \vect(\bm E - \bm G)}_2 \\ 
    &\qquad \geq \rho(\bm H) \norm{\tilde{\bm e}}_2 - \norm{\lambda_1 \bm I_{n^2} - \bm \Lambda_{\bm H}}_2 \norm{\bm Q^{-1}_{\bm H}}_2 \norm{\vect(\bm E - \bm G)}_2 .
\end{align*}
Dividing both sides by $\norm{\tilde{\bm e}}_2$ yields
\begin{align*}
    \frac{\norm{\bm \Lambda_{\bm H} \tilde{\bm e}}_2}{\norm{\tilde{\bm e}}_2} \geq \rho(\bm H) - \frac{\norm{\lambda_1 \bm I_{n^2} - \bm \Lambda_{\bm H}}_2 \norm{\bm Q^{-1}_{\bm H}}_2 \norm{\vect(\bm E - \bm G)}_2}{\norm{\tilde{\bm e}}_2} . \numberthis \label{equ:geqH}
\end{align*}
Since $\norm{\bm E - \bm G}_F = \mathcal{O}(\delta^2)$, (\ref{equ:geqH}) can be rewritten as
\begin{align*}
    \frac{\norm{\bm \Lambda_{\bm H} \tilde{\bm e}}_2}{\norm{\tilde{\bm e}}_2} \geq \rho(\bm H) - \mathcal{O}(\delta^2) . \numberthis \label{equ:geqH1}
\end{align*}
On the other hand, for any $\tilde{\bm e} \in \tilde{\E}_{\bm H}$, we also have 
\begin{align*}
    \frac{\norm{\bm \Lambda_{\bm H} \tilde{\bm e}}_2}{\norm{\tilde{\bm e}}_2} \leq \rho^{\E}(\bm H,\delta) \leq \rho(\bm H) . \numberthis \label{equ:geqH2}
\end{align*}
Combining (\ref{equ:geqH1}) and (\ref{equ:geqH2}) yields $\rho(\bm H) - \rho^{\E}(\bm H,\delta) = \mathcal{O}(\delta)$.

\subsection{Proof of Lemma~\ref{lem:PSD}}
Denote $\bm P_{\bm U} = \bm U \bm U^\topnew$, for any $\bm v \in \R^n$, we can decompose $\bm v$ into two orthogonal component:
\begin{align*}
    \bm v = \bm v_{\bm U} + \bm v_{\perp},
\end{align*}
where $\bm v_{\bm U} = \bm P_{\bm U} \bm v$ and $\bm v_{\perp} = \bm P_{\bm U_\perp} \bm v$. Without loss of generality, assume that $\norm{\bm v}_2 = \norm{\bm v_{\bm U}}_2^2 + \norm{\bm v_{\perp}}_2^2 = 1$. Thus, we have
\begin{align*}
    \bm v^\topnew (\bm M + \bm G) \bm v &= (\bm v_{\bm U} + \bm v_{\perp})^\topnew (\bm M + \bm G) (\bm v_{\bm U} + \bm v_{\perp}) \\
    &= \bm v_{\bm U}^\topnew \bm M \bm v_{\bm U} + \bm v_{\bm U}^\topnew \bm G \bm v_{\bm U} + \bm v_{\bm U}^\topnew \bm G \bm v_{\perp} \\
    &\quad + \bm v_{\perp}^\topnew \bm G \bm v_{\bm U} + \bm v_{\perp}^\topnew \bm G \bm v_{\perp} , \numberthis \label{equ:vM1}
\end{align*}
where the last equation stems from the fact that $\bm M = \P_{\bm U} \bm M \bm P_{\bm U}$ and $\bm P_{\bm U} \bm P_{\bm U_\perp} = \bm 0$. Since $\bm P_{\bm U_\perp} \bm G \bm P_{\bm U_\perp} = \bm 0$, we have
\begin{align*}
    \bm v_{\perp}^\topnew \bm G \bm v_{\perp} = \bm v^\topnew \bm P_{\bm U_\perp} \bm G \bm P_{\bm U_\perp} \bm v = 0 .
\end{align*}
Thus, (\ref{equ:vM1}) is equivalent to
\begin{align} \label{equ:vM2}
    \bm v^\topnew (\bm M + \bm G) \bm v = \bm v_{\bm U}^\topnew \bm M \bm v_{\bm U} + \bm v_{\bm U}^\topnew \bm G \bm v_{\bm U} + 2 \bm v_{\bm U}^\topnew \bm G \bm v_{\perp} .
\end{align}
Now let us lower-bound each term on the RHS of (\ref{equ:vM2}) as follows. First, by the Rayleigh quotient, we have
\begin{align} \label{equ:vM21}
    \bm v_{\bm U}^\topnew \bm M \bm v_{\bm U} \geq \lambda_r \norm{\bm v_{\bm U}}_2^2 ,
\end{align}
and
\begin{align} \label{equ:vM22}
    \bm v_{\bm U}^\topnew \bm G \bm v_{\bm U} \geq \lambda_{\min}(\bm G) \norm{\bm v_{\bm U}}_2^2 \geq - \norm{\bm G}_F \norm{\bm v_{\bm U}}_2^2 .
\end{align}
Next, by Cauchy-Schwarz inequality, 
\begin{align} \label{equ:vM23}
    \bm v_{\bm U}^\topnew \bm G \bm v_{\perp} \geq - \norm{\bm G}_2 \norm{\bm v_{\bm U}}_2 \norm{\bm v_{\perp}}_2 \geq - \norm{\bm G}_F \norm{\bm v_{\bm U}}_2 .
\end{align}
From (\ref{equ:vM21}), (\ref{equ:vM22}), and (\ref{equ:vM23}), we obtain
\begin{align*}
    \bm v^\topnew (\bm M + \bm G) \bm v \geq (\lambda_r - \norm{\bm G}_F) \norm{\bm v_{\bm U}}_2^2 - 2 \norm{\bm G}_F \norm{\bm v_{\bm U}}_2 . \numberthis \label{equ:vMG}
\end{align*}
Note that $\norm{\bm G}_F = \delta$ and the quadratic $g(t) = (\lambda_r - \delta) t^2 - 2 \delta t$ is minimized at 
\begin{align*}
    t_* = \frac{\delta}{\lambda_r - \delta} , \quad g(t_*) = -\frac{\delta^2}{\lambda_r - \delta} .
\end{align*}
Combining this with (\ref{equ:vMG}) yields
\begin{align*}
    \bm v^\topnew (\bm M + \bm G) \bm v \geq -\frac{2}{\lambda_r} \delta^2 ,
\end{align*}
for sufficiently small $\delta$.
Let $\bm F = \bm G + \frac{2}{\lambda_r} \delta^2 \bm I_{n}$. Now we can easily verify that $\norm{\bm F - \bm G}_F = \mathcal{O}(\delta^2)$ and $\bm F \in \F$ .

\subsection{Proof of Lemma~\ref{lem:Pr}}
We shall show that the matrix $\bm E = \P_r(\bm M + \bm F) - \bm M$ belongs to $\E$ and satisfies 
\begin{align} \label{equ:EF}
    \norm{\bm E - \bm F}_F = \mathcal{O}(\delta^2) .
\end{align}
First, since $\bm F \in \F_\delta$, $\bm M + \bm F$ must be PSD. Thus, $\P_r(\bm M + \bm F)$ is a PSD matrix of rank no greater than $r$ and it admits a rank-$r$ factorization $\P_r(\bm M + \bm F) = \bm Z \bm Z^\topnew$, for some $\bm Z \in \R^{n \times r}$. Therefore, by the definition of $\E$,
$$ \bm E = \P_r(\bm M + \bm F) - \bm M = \bm Z \bm Z^\topnew - \bm M \in \E .$$
Next, using (\ref{equ:Pr}), we have
\begin{align*}
    \bm E - \bm F &= \P_r(\bm M + \bm F) - \bm M - \bm F \\
    &= \bm P_{\bm U_\perp} \bm F \bm P_{\bm U_\perp} + \bm{\mathcal{O}}(\norm{\bm F}_F^2) .
\end{align*}
Since $\bm F \in \F_\delta$ implies $\bm P_{\bm U_\perp} \bm F \bm P_{\bm U_\perp} = \bm{\mathcal{O}}(\norm{\bm F}_F^2)$, we conclude that $\bm E - \bm F = \bm{\mathcal{O}}(\norm{\bm F}_F^2)$.



\bibliographystyle{IEEEbib}
\bibliography{refs}

\end{document}